\documentclass{article}[11pt]
\usepackage{amsmath}
  \usepackage{paralist}
  \usepackage{amssymb,amsthm,stmaryrd}
 \usepackage[colorlinks=true]{hyperref}
\hypersetup{urlcolor=blue, citecolor=red}

\newtheoremstyle{drem}
     {3pt}
     {3pt}
     {\rmfamily}
     {}
     {\itshape}
     {:}
     {.5em}
     {}

\newtheorem{teo}{Theorem}[section]
\newtheorem{lem}[teo]{Lemma}
\newtheorem{prop}[teo]{Proposition}
\newtheorem{cor}[teo]{Corollary}
\theoremstyle{drem}
\newtheorem{ex}[teo]{Example}

{\theoremstyle{definition}
\newtheorem{defi}[teo]{{ \textbf{Definition}}}

}

\def\og{\leavevmode\raise.3ex\hbox{$\scriptscriptstyle\langle\!\langle$~}}
\def\fg{\leavevmode\raise.3ex\hbox{~$\!\scriptscriptstyle\,\rangle\!\rangle$}}

\def\vs{\vrule width 0cm height 0.1in depth 0in}
\def\fr#1#2{\frac{\displaystyle\vs #1}{\displaystyle\vs #2}}

\def\bint#1#2{ {\displaystyle\vs \int_{#1}^{#2}} }

\def\ie{\emph{i.e. }}

\def\cf{\emph{cf. }}

\def\eps{\epsilon}

\def\vide{\varnothing}
\renewcommand{\setminus}{\smallsetminus}

\def\un{1 \!\! \mathrm{l}} 
\def\del{\partial \!}

\def\rr{\mathbb{R}}
\def\zz{\mathbb{Z}}
\def\nn{\mathbb{N}}

\def\srl#1{\overline{#1}}
\def\ssl#1{\underline{#1}}

\def\wt#1{\widetilde{#1}}

\def\jo#1{\mathcal{#1}}

\def\limm#1{\textrm{\raisebox{.5ex}{\mbox{$\underset{#1}{\lim}$}}} \:}
\def\lims#1{\textrm{\raisebox{.5ex}{\mbox{$\underset{#1}{\limsup}$}}} \:}
\def\limi#1{\textrm{\raisebox{.5ex}{\mbox{$\underset{#1}{\liminf}$}}} \:}

\def\supp#1{\textrm{\raisebox{.5ex}{\mbox{$\underset{#1}{\sup}$}}} \:}

\newlength{\entrh}   

\def\nr#1{\left\| #1 \right\|}
\def\somme#1#2{\overset{#2}{\underset{#1}{\sum}}}

\def\dd{\mathrm{d} \!}

\def\somme#1#2{\overset{#2}{\underset{#1}{\sum}}}
\def\Somme#1#2{\overset{#2}{\underset{#1}{\displaystyle\vs \sum}}}
\def\produ#1#2{\overset{#2}{\underset{#1}{\prod}}}

\def\union#1#2{\overset{#2}{\underset{#1}{\cup}}}
\def\inter#1#2{\overset{#2}{\underset{#1}{\cap}}}

\def\nr#1{\left\| #1 \right\|}
\def\pnr#1{\| #1 \|}

\def\Bnr#1{\Big\| #1 \Big\|}

\def\abs#1{\left\lvert #1 \right\rvert}

\def\gen#1{\left\langle #1 \right\rangle}
\def\pgen#1{\langle #1 \rangle}

\def\ssi{\Leftrightarrow}
\def\imp{\Rightarrow}
\def\inj{\hookrightarrow}

\def\diam{\mathrm{Diam}\,}
\def\dim{\mathrm{dim}\,}
\def\tr{\mathrm{Tr}\,}
\def\img{\mathrm{Im}\,}
\def\Id{\mathrm{Id}}

\def\homo{\mathrm{Hom}}

\def\IN{\mathrm{int}}
\def\FE{\mathrm{clo}}

\def\wdm{\mathrm{wdim} \,}

\def\mme{\mathrm{Wgc}}

\def\ev{\, e \! v}

\def\dlp{\mathrm{dim}_{\ell^p}}
\def\dl#1{\mathrm{dim}_{\ell^{#1}}}

\title{}
\date{~}

\begin{document}

\centerline{ {\Large A dynamical approach to von Neumann dimension} \footnote[1]{Subject classification, Primary: 37A35, 43A65, 70G60; Secondary: 37C45, 41A46, 43A15, 46E30. }} 
\centerline{\scshape Antoine Gournay \footnote[2]{Email: gournay@mpim-bonn.mpg.de}}
\medskip
{\footnotesize
 \centerline{ Max Planck Institut f\"ur Mathematik}
   \centerline{Vivatsgasse 7}
   \centerline{53111 Bonn, Germany}
} 

\bigskip

\begin{abstract}
Let $\Gamma$ be an amenable group and $V$ be a finite dimensional vector space. Gromov pointed out that the von Neumann dimension of linear subspaces of $\ell^2(\Gamma;V)$ (with respect to $\Gamma$) can be obtained by looking at a growth factor for a dynamical (pseudo-)distance. This dynamical point of view (reminiscent of metric entropy) does not requires a Hilbertian structure. It is used in this article to associate to a $\Gamma$-invariant linear subspaces $Y$ of $\ell^p(\Gamma;V)$ a real positive number $\dlp Y$ (which is the von Neumann dimension when $p=2$). By analogy with von Neumann dimension, the properties of this quantity are explored to conclude that there can be no injective $\Gamma$-equivariant linear map of finite-type from $\ell^p(\Gamma;V) \to \ell^p(\Gamma; V')$ if $\dim V > \dim V'$. A generalization of the Ornstein-Weiss lemma is developed along the way.
\end{abstract}

\section{Introduction}\label{intro}

 Let $\Gamma$ be a discrete group, then it is possible to associate to certain unitary representations a positive real number called von Neumann dimension (see \cite[{\S}1]{Lu} or \cite[{\S}1]{PP}). More precisely, let $f : \Gamma \to X$ be a map. The natural (right) action of $\Gamma$ on spaces of maps means is, in the present text, the action given by $\gamma f(\cdot) = f(\gamma^{-1} \cdot)$. Now, let $H$ be a Hilbert space and consider the space $\ell^2(\Gamma) \otimes H$ where $\Gamma$ acts naturally on the first factor and trivially on the second. Then, the von Neumann dimension is defined for $\Gamma$-invariant subspaces of $\ell^2(\Gamma) \otimes H$.
\par Let $V$ be a \emph{finite} dimensional vector space and $\nr{\cdot}$ a norm (the choice of which will not matter as the dimension is finite).  The subject matter of this article are $\Gamma$-invariant linear subspaces of
\[
\ell^p(\Gamma;V) = \ell^p(\Gamma) \otimes V = \{ f: \Gamma \to V | \somme{\gamma \in \Gamma}{} \nr{ f(\gamma)}^p  \textrm{ is finite}\}
\]
for the natural action of $\Gamma$. Misha Gromov asked the following question (see \linebreak \cite[p.353]{Gro}): are $\ell^p(\Gamma,\rr^n)$ and $\ell^p(\Gamma,\rr^m)$ $\Gamma$-isomorphic if and only if $m=n$?
\par From now on $\Gamma$ will be assumed amenable. There are reasons to exclude non-amenable groups. Indeed, D. Gaboriau pointed out that if a notion of dimension existed in the $\ell^p$ setting (that is, a quantity satisfying properties P1-P10 listed below), then there would be a formula for the Euler characteristic of $\Gamma$ as the alternate sum of the dimensions of $\ell^p$ cohomology spaces. On one hand, torsion-free cocompact lattices in $SO(4,1)$ have positive Euler characteristic. On the other, for $p$ big enough, their $\ell^p$ cohomology vanishes in all degrees but the first (see \linebreak \cite[Theorem 2.1]{PP2}). This would lead to a contradiction.

\par Hence, we are looking for a notion of dimension for such subspaces, which would increase under injective equivariant linear maps. Inspired by an argument of \linebreak \cite[{\S}1.12]{Gro} (and partially answering the question found therein), we shall introduce a quantity $\dlp$ which, when $p=2$, coincides with  definition of von Neumann dimension. This quantity is obtained by a process similar to that of metric entropy or mean dimension, \ie by looking at an asymptotic growth factor. The definition relies \emph{a priori} on an exhaustion of $\Gamma$, but a generalization of the Ornstein-Weiss lemma in section \ref{slow} implies the result is independent of this choice.
\par Though we prove many properties of $\dlp$, important properties are still lacking. Nevertheless, the results obtained in this paper suffice to establish a non existence result for maps of finite type. We recall their construction.
\par Let $D\subset \Gamma$ be a finite set and let $\ssl{g}: V^D \to V'$ be a continuous map. This data enables the definition of a $\Gamma$-equivariant continuous map $g_D$ from $Z \subset \ell^p(\Gamma;V)$ to $\ell^p(\Gamma;V')$ as follows
\[
g_D(z)(\gamma) = \ssl{g}(z(\gamma \delta))_{\delta \in D}.
\]
Remark that what we denote here as $\ell^p(\Gamma;V)$ is more frequently written $\ell^p(\Gamma) \otimes V$.
\begin{teo}\label{t1}
Let $\Gamma$ be an amenable discrete group. Let $V$ and $V'$ be finite dimensional vector spaces. If $f:\ell^p(\Gamma;V) \to \ell^p(\Gamma;V')$ is an injective $\Gamma$-equivariant linear map of finite type then  $\dim V \leq \dim V'$.
\end{teo}
Consequently, if we restrict ourselves to maps of finite type, the question above has a positive answer: there is a $\Gamma$-isomorphism of finite type between $\ell^p(\Gamma;\rr^n)$ and $\ell^p(\Gamma;\rr^m)$ if and only if $m=n$.

\section{Definition and properties of $\dlp$}

Given a positive number $\epsilon$, a notion of dimension up to scale $\epsilon$ for $(X,\tau,\delta)$ a topological space equipped with a pseudo-distance will be needed. Data compression problems turn out to be a good source of inspiration. When one is interested in compression algorithms, it is not only important that the compression map has ``small'' fibers (so that not too much data is lost) but also has an image which is ``small'' in some sense (so that the compression is effective).
\par A slight variant of the one used in \cite{Coo}, \cite{mwi}, \cite{Gro}, \cite{tsu} or \cite{tsu2} shall be employed, namely one that is also defined for pseudo-distances. As such, it will be useful to use a topology $\tau$ that does not come from the pseudo-distance. Please note that the term diameter (denoted $\diam$) will continue to be used even if it is defined using a pseudo-distance (thus a set of diameter $0$ may contain more than one point).
\begin{defi}\label{defwi}
Let $(X,\tau,\delta)$ be a metric space. Call $\wdm_\eps (X,\tau, \delta)$ the smallest integer $k$ such that there exists a continuous (for $\tau$) map $f:X \to K$ where $K$ is a $k$-dimensional polyhedron such that $\forall k \in K, \diam f^{-1}(k) \leq \eps $.
\[
\wdm_\eps (X,\tau,\delta) = \inf_{f:X \inj K } \bigg\{\dim K \bigg|
\begin{array}{c}
f \emph{ is continuous for } \tau \emph{ and } \\
\forall k \in K, \diam f^{-1}(k) < \eps
\end{array}
\bigg\}.
\]
\end{defi}
We will sometimes omit to mention $\tau$ when it is the topology induced by $\delta$.
\begin{defi}
Let $(X,\tau,\delta)$ be a space endowed with a topology $\tau$ and a pseudo-distance $\delta$. Let $\Gamma$ be a countable group which acts on $X$ and let $\{\Omega_i\}$ be an increasing sequence of finite subsets of $\Gamma$. The $\ell^p(\Gamma)$ width growth coefficient of $X$ for the sequence $\{\Omega_i\}$ is
\[
\mme_{\ell^p} (X, \tau, \{\Omega_i\}) = \limm{\eps \to 0} \lims{i \to \infty} \frac{\wdm_\eps (X, \tau, \delta_{\ell^p(\Omega_i)})} {|\Omega_i|} \quad \in [0,+\infty].
\]
\[
\begin{array}{rll}
\textrm{where}  & \delta_{\ell^p(\Omega)}(x,x')= \big(\somme{\gamma \in \Omega}{} \delta(\gamma x, \gamma x')^p \big)^{1/p} & \textrm{for } p<\infty \\
\textrm{and}    & \delta_{\ell^\infty(\Omega)}(x,x') = \supp{\gamma \in \Omega} \delta(\gamma x, \gamma x').
\end{array}
\]
\end{defi}
When $\delta$ is a distance, $\delta_{\ell^\infty(\Omega)}$ is often called the dynamical distance. If furthermore $\tau$ is the topology this distance induces, then this is the (metric) mean dimension (see \cite[{\S}1.5]{Gro} or \cite[{\S}4]{LW}). In the present text, this is an intermediate definition and will only be used in a particular context, namely when $X$ is a subset of $\ell^\infty(\Gamma;V)$. The pseudo-metric will be given by evaluation at the neutral element $e_\Gamma$ of $\Gamma$: $\ev(x,x') = \nr{x(e_\Gamma)-x'(e_\Gamma)}_V$. Lastly, $\tau^*$ will denote the product topology induced from $X \subset V^\Gamma$ (which coincides with the weak-$*$ topology, when defined).
\begin{defi} \label{ddlp}
Let $V$ be a finite-dimensional normed vector space. Let $Y \subset \ell^\infty(\Gamma;V)$ be a subset invariant by the natural action of $\Gamma$, an amenable countable group. Let $\Omega_i$ be a F{\o}lner sequence for $\Gamma$. Then, the $\ell^p$ von Neumann dimension of $Y$ is defined by
\[
\dlp (Y, \{\Omega_i\}) = \supp{r \in \rr_{\geq 0}} \mme_{\ell^p} (B^{Y,p}_r, \ev, \{\Omega_i\})
\]
where $B^{Y,p}_r = Y \cap B^{\ell^p(\Gamma;V)}_r$.
\end{defi}
Note that $B^{Y,p}_r$ is defined by an intersection rather than a projection, as the former are not always easy to define in $\ell^p$. Also, the choice of $\tau^*$ as a topology comes from the fact that it is the weakest topology that is stronger than the topologies induced by $\delta_{\ell^p(\Omega)}$ for any $p$ or $\Omega$.
\par From now on, $Y$ will almost always be a linear subspace. In these cases, one does not need to take the $\sup$ on $r$. Indeed, $\mme_{\ell^p} (B^{Y,p}_r, \ev, \{\Omega_i\})$ does not depend on $r$ (as can be seen using dilation and a change of variable $\eps \mapsto r\eps$).
\par When $Y$ is a $\Gamma$-invariant linear subspace of $\ell^\infty(\Gamma;V)$,
\begin{enumerate}\renewcommand{\labelenumi}{{\normalfont P\arabic{enumi}}}
\item (Independence) $\dlp (Y, \{\Omega_i\})$ is actually independent of the choice of F{\o}lner sequence $\{\Omega_i\}$ (\cf corollary \ref{folind});
\item (Normalization) $\dlp \ell^p(\Gamma;V) =\dim V$ (\cf example \ref{exdimlp});
\item (Invariance) If $p\neq \infty$ and $f:Y_1 \to Y_2$ is an injective $\Gamma$-equivariant linear map of finite type and closed image, then $\dlp Y_1 \leq \dlp Y_2$ (\cf proposition \ref{dlpinv} and corollary \ref{etf});
\item (Completion) If $\srl{Y}$ is the completion of $Y$ in $\ell^p(\Gamma;V)$ for the $\ell^p$ norm, then $\dlp Y =\dlp \srl{Y}$ (\cf proposition \ref{dlpcomp});
\item (Reduction) If $\Gamma_1 \subset \Gamma_2$ is of finite index, and if $Y \subset \ell^p(\Gamma_2;V)$ is seen by restriction as a subspace of $\ell^p(\Gamma_1;V^{[\Gamma_2:\Gamma_1]})$ then $[\Gamma_2:\Gamma_1] \dlp (Y,\Gamma_2) = \dlp (Y,\Gamma_1)$ (\cf proposition \ref{dlpredu}).
\item If $Y \subset \ell^2(\Gamma;V)$, $\dl{2} Y$ coincides with the von Neumann dimension (\cf corollary \ref{wdmVN});
\end{enumerate}
In light of P6, when $p=2$ the following further properties of $\dl{2}$ are listed by Cheeger and Gromov in \cite[{\S}1]{CG}.
\begin{enumerate}\renewcommand{\labelenumi}{{\normalfont P\arabic{enumi}}}\addtocounter{enumi}{6}
\item (Non-triviality) $Y \subset \ell^2$ is trivial if and only if $\dl{2} Y =0$.
\item (Additivity) $\dl{2} Y_1 \oplus Y_2 = \dl{2} Y_1 + \dl{2} Y_2$;
\item (Continuity) If $\{Y_i\}$ is a decreasing sequence of closed linear subspaces then $\dl{2} (\cap Y_i) = \limm{i \to \infty} \dl{2} Y_i$
\item (Reciprocity) If $\Gamma_1 \subset \Gamma_2$ and if $Y_2 \subset \ell^2(\Gamma_2;V)$ is the subspace induced by $Y_1 \subset \ell^2(\Gamma_1;V)$ then $\dl{2} (Y_2,\Gamma_2) = \dl{2} (Y_1,\Gamma_1)$;
\end{enumerate}
\indent Proposition \ref{dlptriv} also establishes P7 for $\dl{1}$. On the other hand, the continuity property (P9) of the von Neumann dimension does not hold if $p=1$ (see example \ref{dl1cex}).
\par For linear subspaces $Y \subset \ell^\infty$ non-triviality (P7) is false, though it might be true for $Y \subset c_0(\Gamma,V)$, the latter being the space of all $x \in \ell^\infty(\Gamma;V)$ tending to $0$ at infinity, \ie $\nr{x}_{\ell^\infty(\Gamma \setminus F_i)} \to 0$ for all exhaustive increasing sequence of (finite) subsets $\{F_i\}$.
\par Finally, the existence of an element of finite support in $Y$ implies P7. By using a similar but less convenient definition of $\dl{2}$, the author is also aware of a proof of P7, P8 and P10 (when the index is finite) for $p=2$ without using P6 and the previously known properties of von Neumann dimension.
\par Though these properties are stated for $\Gamma$-invariant linear subspaces, some remain true for more general subsets $Y$: P1 and P5 hold for any $\Gamma$-invariant subset, and P4 is also true when $Y$ is not $\Gamma$-invariant.
\par These properties are sufficient to offer a partial answer to the question discussed at the beginning of the present article.
\begin{proof}[Proof of theorem \ref{t1}]
It is but a simple consequence of P2 and P3.
\end{proof}
Being crucial to the proof above and less technical, we shall begin by proving properties P2-P5. Section \ref{fpdlp} then discusses P7 and P9 for $p=1$ or $\infty$. The proof of P1 requires some quite technical lemmas on amenable groups and is thus relegated to section \ref{slow}. As for P6, it relies mostly on a result of Gromov and is discussed in appendix A.

\section{Proof of properties P2-P5} \label{sdlprop}

\hspace*{0.75\parindent} Before the properties of $\dlp$ can be established, the basic properties of $\wdm_\eps$ must be mentioned.

\subsection{Properties of $\wdm_\eps$}
Most of the content of this subsection may be found in \cite[{\S}4.5]{Coo}, \cite[{\S}3]{CK}, \cite[Proposition 2.1]{mwi} and \cite[{\S}1.1]{Gro}.
\begin{prop}\label{propbase1} Let $X$ be space endowed with a topology $\tau$ and a pseudo-distance $\delta$.
  \begin{enumerate}\renewcommand{\labelenumi}{{\normalfont \alph{enumi}.}}
    \item The function $\eps \mapsto \wdm_\eps (X,\tau,\delta)$ is non-increasing.
    \item Suppose $\delta$ is a distance and $\tau$ the topology it induces. Let $\dim X$ be the covering dimension of $X$, then $\wdm_\eps (X,\tau,\delta) \leq \dim X$.
    \item $\wdm_\eps (X,\tau,\delta) =0 \ssi \eps \geq \diam X$.
  \end{enumerate}
\end{prop}
Except for b, the proof of these properties are simple. Before moving on, let us recall two (fundamental) examples.
\begin{ex}\label{basm}
Let $X$ be a normed vector space with the distance $\delta(x,x') = \nr{x-x'}$ and $\tau$ the norm topology. Let $A = B^X_1$ be its unit ball. Then (see \cite[{\S}1.1B]{Gro} or, for more details, \cite[Lemma 2.5]{mwi} or \cite[Appendix]{tsu2}) $\wdm_\eps (A,\tau,\delta) = \dim X$ if $\eps < 1$ and  (by considering the map which sends all of $A$ to one point) $\wdm_\eps (A,\tau,\delta) = 0$ if $\eps \geq 2$.
\end{ex}
The second example comes from a question which arises naturally in the context of compressed sensing, namely we look at a ball for some norm but we endow with a metric coming from another norm.
\begin{ex}\label{blplq}
Let $\ell^p(n)$ denote $\rr^n$ with its $\ell^p$ norm. Then one can try to compute the $\wdm$ of $B^{\ell^q(n)}_1 \subset \ell^p(n)$. (In compression theory, it is frequent to consider a ball for some metric endowed with a different metric; see \cite{Dono}.) When $q \geq p$ then the behaviour is essentially as in the previous example. However, if $q<p$ then one finds that, for $1 \leq k < n$,
\[
\wdm_\eps (B_1^{\ell^p(n)},\ell^q)  \left\{
  \begin{array}{llrcl}
    = 0     & \textrm{if} &                2 \leq & \eps, & \\
    \leq k  & \textrm{if} & 2(k+1)^{1/q-1/p} \leq & \eps, & \\
    \geq k  & \textrm{if} &                       & \eps  & < k^{1/q-1/p}, \\
    = n     & \textrm{if} &                       & \eps  &< n^{1/q-1/p}.
  \end{array} \right.
\]
We briefly mention how to obtain these. The first line is a consequence of \ref{propbase1}.c. The second is found by using an explicit map described in \cite[proposition 1.3]{mwi} and \cite{tsu}. This maps takes a vector, keeps only the $k$ biggest coordinates (in absolute value), then from these $k$ coordinates take the smallest and substract (or add, so as to reduce in absolute value) it to the others. Finally, the third line comes from the presence of an $\ell^q$ ball of dimension $k$ and radius $k^{1/q-1/p}$ in $\ell^p(n)$. The fourth line is also obtained using this argument (for $n$) together with proposition \ref{propbase1}.b.
\end{ex}
This second set of properties are crucial to what follows.
\begin{prop}\label{propbase}
For $i=1,2$, let $X_i$ be spaces endowed with topologies $\tau_i$ and pseudo-metric $\delta_i$.
  \begin{enumerate}\renewcommand{\labelenumi}{{\normalfont \alph{enumi}.}}
    \item Let $f: X_1 \to X_2$ be a continuous map such that $\delta_1(x,x') \leq C \delta_2\big(f(x),f(x')\big)$ where $C \in ]0,\infty[$. Then $\wdm_\eps (X_1,\tau_1,\delta_1) \leq \wdm_{\eps/C} (X_2,\tau_2,\delta_2)$.
    \item A dilation has the expected effect, \ie let $f:X_1 \to X_2$ be a homeomorphism such that  $\delta_1(x,x') = C \delta_2\big(f(x),f(x')\big)$. Then $\wdm_\eps X_1 = \wdm_{\eps/C} X_2$.
    \item Let $X:=X_1 \times_q X_2$ be the space $X_1 \times X_2$ endowed with the product topology and the pseudo-metric $\delta:= \delta_1 \times_q \delta_2$ given by $\delta(x,x')^q = \delta_1(x,x')^q + \delta_2(x,x')^q$, when $q \in [1,\infty[$, and $\delta(x,x') = \max \big( \delta_1(x'x'), \delta_2(x,x') \big)$ when $q=\infty$. Then $\wdm_{2^{1/q}\eps} X \leq \wdm_{\eps} X_1 + \wdm_{\eps} X_2$.
  \end{enumerate}
\end{prop}
The proofs can be found in \cite[{\S}4.5]{Coo}, \cite[Lemma 3.2]{CK} or \cite[Proposition 2.1]{mwi}. For example, the third is obtained by looking at the size of the fibers of the map $f = f_1 \oplus f_2$, where $f_i:X_i \to K_i$ satisfy the conditions of definition \ref{defwi} and $\dim K_i = \wdm_\eps (X_i, \tau_i, \delta_i)$.
\par A useful way of stating \ref{propbase}.a is that a continuous map that does not reduce distances will not make $\wdm_\eps$ smaller.

\subsection{Properties of $\dlp$}
Let us begin by two basic examples.
\begin{ex}\label{exdimlp}
\emph{ If $1\leq q<p \leq \infty$, and $Y=B^{\ell^q(\Gamma;\rr)}_1$ then \mbox{$\dlp (Y, \{\Omega_i\}) =0$} }(independently of the choice of sequence $\{\Omega_i\}$). Indeed, $B^{Y,p}_r \subset B^{\ell^q(\Gamma;\rr)}_{r}$, and $\wdm_\eps (B^{Y,p}_r, \ev_{\ell^p(\Omega_i)}) = \wdm_\eps (B^{\ell^q(n_i)}_{r'},\ell^p)$ where $n_i = |\Omega_i|$. However using example \ref{blplq} (and dilations to get back to a unit ball, see proposition \ref{propbase1}.b), \linebreak $\wdm_\eps (B^{\ell^q(n_i)}_{r'},\ell^p)$ is, for fixed $\eps$, bounded above and below by two functions that do not depend on $n_i$. Thus,
\[
\lims{i\to \infty} \fr{\wdm_\eps (B^{Y,p}_r, \tau^*, \ev_{\Omega_i}  )}{|\Omega_i|} \leq \lims{i\to \infty} \fr{\wdm_\eps (B^{\ell^q(n_i)},\ell^p)}{|\Omega_i|} = 0.
\]
\end{ex}
\begin{ex}\label{exdimlp2}
\emph{By direct computation, we now show that }
\[
\dlp \ell^q(\Gamma;V) = \dim V.
\]
 For $q\in [1,\infty]$, let $Y' = \ell^q(\Gamma;V)$. Then $(B^{Y',p}_1, \tau^*, \ev_{\ell^p(\Omega)})$ is  ``isometric'' (for the pseudo-distance) to $(B^{\ell^p(\Omega;V)}_1, \ev_{\ell^p(\Omega)})$. Indeed, the restriction map to $\Omega$ has a kernel of ``diameter'' $0$, so property \ref{propbase}.a applies with $C=1$. On the other hand, inclusion of $\ell^p(\Omega,V)$ in $\ell^p(\Gamma,V)$ (by extending the functions by $0$) is also a linear map and property \ref{propbase}.a holds again with $C=1$. Consequently, $(B^{Y',p}_1, \tau^*, \ev_{\ell^p(\Omega)})$ will have the same $\wdm_\eps$ as $(B^{\ell^p(\Omega;V)}_1, \ev_{\ell^p(\Omega)})$, $\forall \eps$. This later being a ball with its proper metric, if $\eps<1$ its $\wdm_\eps$ will be the dimension of the space, $|\Omega| \dim V$.
\end{ex}
\par In what follows the total vector space will be $Y \subset \ell^p(\Gamma;V)$. So  $(B^{Y,p}_r, \tau^*, \ev_{\ell^p(\Omega)})$ is the set $B^{Y,p}_r \subset \ell^p(\Gamma;V)$ with the pseudo-norm $\ev_{\ell^p(\Omega)}$ and the topology induced from the product topology. We stress that $B^{Y,p}_r$ is not the ball for the pseudo-norm $\ev_{\ell^p(\Omega)}$; it is the intersection of $Y$ with the ball of radius $r$ in $\ell^p(\Gamma)$ (endowed with its actual norm). The next property is a corollary of a generalization of the Ornstein-Weiss lemma described in section \ref{slow}.
Even if proposition \ref{folind} is a very important property, weaker version can be sufficient for some of our needs. Indeed, the following simple lemma is actually all that we need to show that $\dlp$ is preserved under $\Gamma$-equivariant maps of finite type.
\begin{lem}\label{decsu}
Let $Y$ be as above, and let $\{ \Omega_i \} $ and $\{ \Omega'_i \} $ be such that
\[
\limm{i \to \infty} \fr{|\Omega_i \cup \Omega_i' \setminus \Omega_i \cap \Omega_i'|}{|\Omega_i \cup \Omega_i'|} = 0,
\]
then $\mme_{\ell^p} (B^{Y,p}_1, \ev, \{\Omega_i\}) = \mme_{\ell^p} (B^{Y,p}_1, \ev, \{\Omega'_i\})$
\end{lem}
\begin{proof}
It suffices to note that, when $\Omega \subset \Omega'$,
\[
\begin{array}{rcl}
\fr{\wdm_\eps (B^Y_1, \ev_{\ell^p(\Omega)})}{|\Omega|} \fr{|\Omega|}{|\Omega'|}
   &\leq& \fr{\wdm_\eps (B^Y_1, \ev_{\ell^p(\Omega')})}{|\Omega'|} \\
   &\leq& \fr{\wdm_\eps (B^Y_1, \ev_{\ell^p(\Omega)})}{|\Omega|}\fr{|\Omega|}{|\Omega'|}
            + \dim V \fr{|\Omega' \setminus \Omega|}{|\Omega'|}.
\end{array}
\]
Furthermore, $\fr{|\Omega|}{|\Omega'|} = 1- \fr{|\Omega' \setminus \Omega|}{|\Omega'|}$. Thus, computing $\mme$ with respect to the sequences $\{\Omega_i \cap \Omega_i'\}$, $\{\Omega_i\}$ or $\{\Omega_i'\}$ will yield the same result as a computation made using $\{ \Omega_i \cup \Omega_i'\}$.
\end{proof}
\begin{prop}\label{dlpinv}
Let $Y \subset \ell^\infty(\Gamma;V)$ and $Y'\subset \ell^\infty(\Gamma;V')$ be $\Gamma$-invariant linear subspaces. Let $f: Y \to Y'$ be a $\Gamma$-equivariant map continuous for $\tau^*$ and such that there exists a real $c_f \in \rr_{>0}$ and a finite subset $D_f \subset \Gamma$ satisfying \linebreak $\ev(x,y) \leq c_f \ev_{\ell^p(D_f)}(f(x),f(y))$ then
\[
\dlp (Y,\{\Omega_i\}) \leq \dlp (Y',\{\Omega_i\})
\]
\end{prop}
\begin{proof}
The case $p=\infty$ is simpler, we shall only describe the case $p<\infty$. Here $B^{Y',p}_r = Y' \cap B^{\ell^p(\Gamma;V')}_r$. On one hand, since $f$ is continuous for $\tau^*$ (the product topology or the weak-$^*$ topology), $\exists r_f \in \rr_{>0}$ such that $f(B^Y_1) \subset B^{Y',p}_{r_f}$. Indeed, since the image is weakly-$^*$ compact (in particular, weakly-$^*$ bounded) it is bounded (\cf \cite[theorem 3.18]{Ru}). On the other hand, the assumption satisfied by $f$ on distances propagates by equivariance to different evaluations:
\[
\begin{array}{rl}
\ev(\gamma x, \gamma y) \leq c_f \ev_{\ell^p(D_f)}(f(\gamma x),f(\gamma y))
  =& \hspace*{-1.5ex} c_f \ev_{\ell^p(D_f)}(\gamma f(x), \gamma f(y)) \\
  =& \hspace*{-1.5ex} c_f \ev_{\ell^p(D_f \gamma)}(f(x), f(y)).
\end{array}
\]
This implies that $\ev_{\ell^p(\Omega)}(x,y) \leq c_f |D_f| \ev_{\ell^p(\Omega D_f)}(f(x),f(y))$ and, incidentally, that $f$ is injective. Lastly, since the image of the ball (of radius $1$) is contained in a ball (of radius $r_f$)
\[
\begin{array}{rl}
\wdm_\eps (B^{Y,p}_1, \ev_{\ell^p(\Omega_i)})
   \leq& \hspace*{-1.5ex} \wdm_{\eps/c_f |D_f| } (B^{Y',p}_{r_f},\ev_{\ell^p(D_f \Omega_i)}) \\
   \leq& \hspace*{-1.5ex} \wdm_{\eps /c_f |D_f| r_f} (B^{Y',p}_1, \ev_{\ell^p(\Omega_i D_f)}).
\end{array}
\]
The first inequality comes from \ref{propbase}.a. Dividing by $|D_f \Omega_i| = \frac{|D_f \Omega_i|}{|\Omega_i|} |\Omega_i|$ and passing to the limit yields that
\[
\mme_{\ell^p} (B^{Y,p}_1,\ev,\{\Omega_i\}) \limm{i \to \infty} \frac{|D_f \Omega_i|}{|\Omega_i|} \leq \mme_{\ell^p} (B^{Y',p}_1,\ev,\{\Omega_i D_f \}).
\]
Since $\{\Omega_i \}$ is a F{\o}lner sequence, the limit on the left-hand side is $1$. Furthermore, the hypothesis of lemma \ref{decsu} are satisfied; the right-hand term is nothing else than $\dlp (Y', \{\Omega_i\})$.
\end{proof}
From now on, we will drop the explicit reference to the F{\o}lner sequence.
\par Since the assumptions of the previous proposition are quite abstract, it is good to check that they hold in certain categories of maps. The main constraint is the existence of $c_f$ and $D_f$.
Let $f$ be a map to which proposition \ref{dlpinv} applies. Let $f^{-1}: Y' = \img f \to Y$ the inverse of $f$ on its image, then the condition
\[
\ev(x,y) \leq c_f \ev_{\ell^p(D_f)}(f(x),f(y))
\]
can be read as a condition on the modulus of continuity of $f^{-1}$. More precisely, \linebreak $f^{-1}:(Y',\ev_{\ell^p(D_f)}) \to (Y,\ev)$ must be continuous with a linear modulus of continuity, \ie $f^{-1}$ must be Lipschitz. If the function $f^{-1}$ is continuous for the product topology, weakening the topology on its image is evidently not restrictive. Things are not so direct on the domain.
\par For $\Omega \subset \Gamma$, denote by $R_\Omega : V^\Gamma \to V^\Omega$ the restriction of functions to the domain $\Omega$. Let $U \subset (Y,\ev)$ be an open set; then if $Y$ is seen as a subset of $V^\Gamma$, $R_{\{e_\Gamma\}} U$ is an open set on the factor $R_{\{e_\Gamma\}} Y$, and all of $Y$ on the other factors: $R_{\Gamma \setminus \{e_\Gamma\}} U = R_{\Gamma \setminus \{e_\Gamma\}} Y$. It is then possible that on a finite number of factors of $Y' \subset V'^\Gamma$ (the required set $D_f$) $f(U)$ will not be all the image of $f$: $R_D f(U) \neq R_D Y'$. For example, for $F\subset \Gamma$ a finite subset and $f = f_F$ of finite type, the condition is that $f:(Y,\ev) \to (Y,\ev_{\ell^p(F)})$ be open (on its image) and of Lipschitz inverse. Remember that the condition on the distances in proposition \ref{dlpinv} and $\Gamma$-equivariance imply injectivity of the map. Here is the major application of proposition \ref{dlpinv}.
\begin{cor}\label{etf}
Let $Y \subset c_0(\Gamma;V)$ and $Y'\subset c_0(\Gamma;V')$ be $\Gamma$-invariant linear subspaces. Let $f: Y \to Y'$ be a $\Gamma$-equivariant injective linear map of finite type and closed image. Then
\[
\dlp Y \leq \dlp Y'.
\]
\end{cor}
\begin{proof}
Let $F \subset \Gamma$ be a finite subset which can be used to define $f$ as a map of finite type, \ie $f=f_F$. If $f$ is a linear map of closed image, injectivity of $f$ implies that it is open on its image (Banach-Schauder theorem or open mapping theorem) for the norm topologies. This remains true for the topology of $\ev$ on the domain and $\ev_\ell^p(F)$ on the image as the first is weaker (its open sets are described above) and the image of open sets is of the form $R_F^{-1} U'$ for $U' \subset V^F$. So $f:(Y,\ev) \to (Y,\ev_{\ell^p(F)})$ is open (on its image).
\par Next, write the $\Gamma$-equivariant linear map of finite type $f$ as
\[
x \mapsto f(x) \textrm{ such that } f(x)(\gamma) = \somme{\gamma' \in F}{} a_{\gamma'} \big( x(\gamma' \gamma) \big),
\]
where $a_{\gamma'} \in \homo(V,V')$. Since it is injective and the finitely supported functions are dense in $Y'$ ($\ell^\infty$ being excluded), it possesses a $\Gamma$-equivariant linear inverse (on its image) $f^{-1} = g$ of the form:
\[
x \mapsto g(x) \textrm{ such that } g(x)(\gamma) = \somme{\gamma' \in G}{} b_{\gamma'} \big( x(\gamma' \gamma) \big),
\]
where $b_{\gamma'} \in \homo(V',V)$ and $G \subset \Gamma$ might not be finite. 
\par Then proposition \ref{dlpinv} can be invoked with $D_f=F$, and $c_f$ the Lipschitz constant of $g:(Y',\ev_{F^{-1}}) \to (Y,\ev)$. Thus $c_f \leq \pnr{\oplus_{\gamma \in F^{-1} \cap G} b_\gamma}$.
\end{proof}
\begin{prop}\label{dlpcomp}
  Let $Y \subset \ell^\infty(\Gamma;V)$ be an open linear subspace and let $\srl{Y}$ be its completion in $\ell^p(\Gamma;V)$, then $\dlp Y = \dlp \srl{Y}$.
\end{prop}
\begin{proof}
  The argument is identical to that of example \ref{exdimlp}: when restricted to a finite $\Omega \subset \Gamma$, these two spaces cannot be distinguished (being of finite dimension they are closed). In other words, there exists a linear map, given by the restriction $R_\Omega$, and whose kernel is in the ``ball'' of radius $0$:
\[
R_\Omega: (B^{\srl{Y}}_1, \ev_{\ell^p(\Omega)})  \to  (R_\Omega B^{Y}_1, \ev_{\ell^p(\Omega)}).
\]
Thus, $\forall \eps \in [0,1]$, $\wdm_\eps(B^{\srl{Y},p}_1, \ev_{\ell^p(\Omega)}) \leq \wdm_\eps(R_\Omega B^{Y,p}_1, \ev_{\ell^p(\Omega)})$.
\par On the other hand, let $s: R_\Omega B^{Y,p}_1 \to B^{Y,p}_1$ such that $R_\Omega \circ s = \Id$ be determined by an inverse of $R_\Omega Y \to Y$, then $s$ is a linear map which increases distances. Consequently, $\wdm_\eps(R_\Omega B^{Y,p}_1, \ev_{\ell^p(\Omega)}) \leq \wdm_\eps(B^{Y,p}_1, \ev_{\ell^p(\Omega)})$. Finally, by inclusion $Y \subset \srl{Y}$, we have $\wdm_\eps(B^{Y,p}_1, \ev_{\ell^p(\Omega)}) \leq  \wdm_\eps(B^{\srl{Y},p}_1, \ev_{\ell^p(\Omega)})$.
\end{proof}
\par If $[\Gamma_2:\Gamma_1] = |G| <\infty$, a set $Y \subset \ell^p(\Gamma_2;V)$ is also a set of $\ell^p(\Gamma_1;V^G)$. Indeed, to $y \in \ell^p(\Gamma_2;V)$ one can associate $i(y)$ where $i(y)(\gamma) = (y(\gamma g))_{g \in G} \in V^G$. This operation behaves nicely with $\dlp$.
\begin{prop} \label{dlpredu}
Let $\Gamma_1 \subset \Gamma_2$ be amenable groups and $G =\Gamma_2 / \Gamma_1$ where $|G|$ is finite, if $Y \subset \ell^p(\Gamma_2;V)$ is seen by restriction as a linear subspace of $\ell^p(\Gamma_1;V^{G})$ then $|G| \dlp (Y,\Gamma_2) = \dlp (Y,\Gamma_1)$.
\end{prop}
\begin{proof}
Let $\{\Omega^{(1)}_i \}$ be a F{\o}lner sequence for $\Gamma_1$ and let $\{\Omega^{(2)}_i \} = \{\Omega^{(1)}_i G \}$ be the corresponding F{\o}lner sequence in $\Gamma_2$. It is then sufficient to see that $(B^{Y_2}_1 , \ev_{\Omega_i^{(2)}})$ is by construction isometric to $(B^{Y_1}_1 , \ev_{\Omega_i^{(1)}})$.
\end{proof}
\indent Let us mention a typical problem when one deals with $\ell^p$ spaces, for $p \neq 2$, that is the existence of linear subspaces which are not the image of projection (\cf \cite{Mur} and \cite{Sob}). A characterization of  subspaces of $\ell^p$ possessing a projection of norm $1$ can be found in \cite[I.{\S}2]{LT}. We shall briefly discuss the case where $Y \subset \ell^p(\Gamma;\rr)$ is a $\Gamma$-invariant linear subspace and there exists a $\Gamma$-equivariant bounded linear map, $P_Y$ (which is not necessarily a projection). Then let $y=P_Y \delta_{e_\Gamma}$ where $\delta_{e_\Gamma}$ is the Dirac mass at ${e_\Gamma} \in \Gamma$, and let $q \leq p$ be such that $y \in \ell^q(\Gamma;\rr)$. For a $x \in \ell^p(\Gamma;\rr)$, write $x = \sum k_\gamma \delta_\gamma$. By linearity and $\Gamma$-equivariance of $P_Y$,
\[
P_Y x = P_Y \somme{\gamma \in \Gamma}{} k_\gamma \delta_\gamma = \somme{\gamma \in \Gamma}{} k_\gamma P_Y(\gamma \delta_{e_\Gamma}) = \somme{\gamma \in \Gamma}{} k_\gamma \gamma y
\]
Choosing $k_\gamma = |y(\gamma^{-1})|^{\frac{q}{p} -1}y(\gamma^{-1})$, the evaluation $(P_Y x)({e_\Gamma}) = \sum |y(\gamma^{-1})|^{\frac{q}{p} + 1}$ must be finite. This forces $\frac{q}{p} + 1 \leq q$, in other words $q \leq p'$ (where $p'$ is the conjugate exponent to $p$). When $p>2$, the existence of such a map means has the (restrictive) consequence that there exists in $Y$ an element which is also in $\ell^{p'}(\Gamma;\rr)$.

\section{Further properties in special cases}\label{fpdlp}

\hspace*{0.75\parindent} We now discuss property P7, that is if $Y$ is non-trivial then $\dlp Y$ is positive. This question is difficult as an intuitive proof only works for $p=1$. Before we move to this proof, let us argue why the three following assumptions seem necessary for it to hold: $Y$ must be a linear subspace, $Y$ must be $\Gamma$-invariant, and $Y$ must be contained in $\ell^p(\Gamma;V)$ for finite $p$ or in $c_0(\Gamma;V)$ if $p=\infty$. Here are some cases of non-trivial $Y$ for which one of the assumptions does not hold and where $\dlp$ is $0$.
\par First, suppose $Y$ is not a linear subspace. In example \ref{exdimlp} the $\ell^q$ balls where $q < p$ are shown to have their $\dlp$ equal to $0$. Alternatively, one could also take $Y$ to be the subset of $\ell^\infty(\Gamma;V)$ given by function with support of cardinality less than $k$ (for a fixed $k \in \zz_{>0}$).
\par Second, if $Y$ is a linear subspace of $\ell^\infty(\Gamma;V)$ but is not $\Gamma$-invariant, it could be of finite dimension, and consequently $\dlp$ will be trivial.
\par Last, when $p$ is finite, the existence of a $y \in Y$ whose $\ell^p$ norm is finite is only guaranteed if $Y \subset \ell^p$. Otherwise, it could happen that $Y \cap B^{\ell^p(\Gamma;V)}_r = \{0\}, \forall r$. On the other hand, if $p=\infty$, take $Y \subset \ell^\infty(\Gamma;V)$ the ($\Gamma$-invariant) line generated by a constant function $y$ (\ie such that $\exists v \in V, \forall \gamma \in \Gamma, y(\gamma) = v$). $Y$ is $1$-dimensional, and consequently $\dl{\infty} Y=0$. But $Y$ is not trivial. However, the question for a $\Gamma$-invariant linear subspace $Y \subset c_0(\Gamma;V)$ remains interesting.
\par Fortunately, in the $\ell^1$ case things can be proved without difficulties. As noted before this method does not extend to $p>1$.
\begin{prop}\label{dlptriv}
Let $Y \subset \ell^1(\Gamma;V)$ be a $\Gamma$-invariant linear subspace, then \linebreak $\dl{1} Y = 0$ if and only if $Y$ is trivial.
\end{prop}
\begin{proof}
This proof requires some results on amenable groups; these can be found in section \ref{slow}. If one wants, it is possible to think of $\Gamma$ as $\zz^n$ and take finite sets to be rectangles.
\par If $Y$ is trivial then $\dl{1} Y$ is obviously $0$. Otherwise, let $0\neq y \in Y$ and renormalize it so that $\nr{y}_{\ell^1(\Gamma)} = 1$. For all $\eps \in ]0,1/2[$, $\exists F \subset \Gamma$ finite (which depends on $y$ and $\eps$) such that $\nr{y}_{\ell^1(F)} > 1-\eps$ (and consequently $\nr{y}_{\ell^1(\Gamma \setminus F)} \leq \eps$). Then let $\wt{y}$ be identical to $y$ on $F$ and $0$ elsewhere.
\par For $i$ sufficiently big, $\Omega_i$ contains a non-empty $\rho$-quasi-tiling by $F$ (see definition \ref{dqt}), since $F\subset \Omega_i$ and $\alpha(\Omega_i;F)$ tends to $0$. Applying lemma \ref{epac} to find translates of $F$ which are $\rho$-disjoint, where $\rho= 1/2|F|$, we obtain a quasi-tiling whose elements are actually disjoint since $\rho < |F|^{-1}$, and the number of such translates is at least $(1-\alpha(\Omega_i;F))|\Omega_i|/2|F|$.
\par Let $\gamma_j$ for $j \in J_i \subset \zz_{>0}$ be the elements by which the sets $F$ are translated for a $\rho$-quasi-tiling of $\Omega_i$ (since the $\Omega_i$ form an increasing sequence and that lemma \ref{epac} applies to all maximal $\rho$-quasi-tiling, it can be assumed that the $J_i$ are increasing). Let $V_i = \gen{\gamma_j y| j \in J_i}$ be the linear subspace generated by the corresponding translates of $y$. Trivially $B_1^{V_i} \subset B^Y_1$, and  we will construct a map from a ball to $B^{V_i}_1$. Let
\[
\begin{array}{rcll}
  \pi:& \ell^1(J_i;\rr)   &\to    & V_i\\
      &  (a_j)_{j \in J_i}&\mapsto&  \somme{j \in J_i}{} a_j \gamma_j y
\end{array}
\textrm{ and }
\begin{array}{rcll}
\wt{\pi}:& \ell^1(J_i;\rr)   &\to    & V_i\\
         &  (a_j)_{j \in J_i}&\mapsto&  \somme{j \in J_i}{} a_j \gamma_j \wt{y}
\end{array}
\]
With these notations,
\[
\begin{array}{rll}
\nr{\wt{\pi}(a)}_{\ell^1(\Gamma)}
    &= \somme{k \in J_i}{} \Bnr{\somme{j \in J_i}{} a_j \gamma_j \wt{y} }_{\ell^1(\gamma_k F)}
    &= \somme{k \in J_i}{} \Bnr{a_k \gamma_k \wt{y} }_{\ell^1(\gamma_k F)} \\
    &= \somme{k \in J_i}{} |a_k| \nr{\wt{y}}_{\ell^1(F)}
    &= \nr{y}_{\ell^1(F)} \somme{k \in J_i}{} |a_k| .
\end{array}
\]
On the other hand,
\[
\begin{array}{rll}
\nr{\wt{\pi}(a) - \pi(a)}_{\ell^1(\Gamma)}
   &= \Bnr{ \somme{j \in J_i}{} a_j \gamma_j(\wt{y}-y) }_{\ell^1(\Gamma)}
   &= \somme{\gamma \in \Gamma}{} \Big| \somme{j \in J_i}{} a_j \gamma_j(\wt{y}-y) \Big| \\
   &\leq \somme{\gamma \in \Gamma}{} \somme{j \in J_i}{} |a_j \gamma_j(\wt{y}(\gamma)-y(\gamma))|
   &= \somme{\gamma \in \Gamma}{} \somme{j \in J_i}{} |a_j| |\wt{y}(\gamma)-y(\gamma)|  \\
   &= \somme{j \in J_i}{} |a_j| \Big( \somme{\gamma \in \Gamma}{} |\gamma_j(\wt{y}(\gamma)-y(\gamma))| \Big)
   &= \nr{y}_{\ell^1(\Gamma \setminus F)} \somme{j \in J_i}{} |a_j|.\\
\end{array}
\]
In short,
\[
\nr{\wt{\pi}(a)}_{\ell^1(\Gamma)} = \nr{y}_{\ell^1(F)} \nr{a}_{\ell^1(J_i)} \textrm{ and } \nr{\wt{\pi}(a) - \pi(a)}_{\ell^1(\Gamma)} \leq \nr{y}_{\ell^1(\Gamma \setminus F)} \nr{a}_{\ell^1(J_i)}.
\]
Since $1-\eps \leq \nr{y}_{\ell^1(F)} \leq 1$ and $\nr{y}_{\ell^1(\Gamma \setminus F)} \leq \eps $, we conclude that
\[
(1-2\eps) \nr{a}_{\ell^1(J_i)} \leq \nr{\pi(a)}_{\ell^1(\Gamma;V)}\leq ( 1+\eps) \nr{a}_{\ell^1(J_i)}
\]
This means that $(B^Y_1,\ev_{\ell^1(\Omega_i)})$ contains, with a controlled distortion, a $\ell^1$ ball (with its $\ell^1$ metric) of radius $1$ and of dimension $\frac{1}{2|F|} (1-\alpha(\Omega_i;F))|\Omega_i|$, whence
\[
\begin{array}{rl}
\dl{1}Y
   &= \limm{\eps' \to 0} \lims{i \to \infty} \frac{\wdm_{\eps'} (B^Y_1,\ev_{\ell^1(\Omega_i)}) }{|\Omega_i|} \\
   &\geq \limm{\eps' \to 0} \lims{i \to \infty} \frac{1}{2|F|} (1-\alpha(\Omega_i;F)) = \frac{1}{2|F|}.
\end{array}
\]
As required $\dl{1}Y > 0$.
\end{proof}
\par This result can be extended to $p>1$ in the special case that $Y \subset \ell^p(\Gamma;V)$ contains an element in $\ell^1$ (in particular, an element of finite support). By P6, positivity (P7) is also true for $p=2$. Positivity means that if one looks at a \linebreak $p$-summable two-sided sequence $y \in \ell^p(\zz;\rr)$, there are subspaces of the space generated by $y$ and sequences obtained by shifting $y$ up to $n$ times left or right of dimension proportional to $n$ and so that no element decreases too much in norm when restricted to $[-n,n]$. The above result is a simple consequence that this is true for $p=1$, $p=2$ is also true albeit not so simple, and one is then lead to ask if this can be true for other values of $p\neq \infty$ or in $c_0$.
\par Even if we cannot show continuity, the following example is worthy of interest. The sequence of vector subspaces discussed there will not satisfy the continuity property (P9). This is quite unfortunate, as $\ell^1$ is among the few cases where positivity can be shown.
\begin{ex}\label{dl1cex}
\emph{We exhibit a decreasing sequence of closed linear subspaces of $\ell^1(\zz,\rr)$, $\{Y_i\}$, such that}
\[
1 = \limm{i \to \infty} \dl{1} Y_i \neq \dl{1} \inter{i \to \infty}{} Y_i =0.
\]
Define $\forall k \in \zz_{>0}, \pi_k: \ell^1(\zz;\rr) \to \ell^\infty(\zz/k\zz; \rr)$ in the following way: for $n \in \zz/k\zz$
\[
\pi_k(x)(n) = \somme{i \equiv n \textrm{ mod }k}{} x(i).
\]
Continuous linear maps between Banach spaces have a closed kernel (for $\tau$, the norm topology in $\ell^1$), thus $Y_j= \inter{k=1}{j} \ker \pi_k$ is a decreasing sequence of closed sets (for $\tau$). To compute $\dl{1}$, choose the F{\o}lner sequence $\Omega_i = [-i,i] \cap  \zz$. For a $N \in \nn$, let $y_N \in Y_1$ be such that $y_N(0)=1/2$, $y_N(N)=-1/2$ and which is zero elsewhere. Let $N_j = \textrm{lcm}(1,2,\ldots,j)$. For all $j$, $y_{N_j} \in B^{Y_j}_1$. These elements give a map $(B^{\ell^1(\zz;\rr)}_{1/2}, \ev_{\ell^1(\Omega)})$ to $(B^{Y_j}_1, \ev_{\ell^1(\Omega)})$ which possesses fibers of ``diameter'' $0$. They are defined as follows, $y \in B^{\ell^1(\zz;\rr)}_{1/2}$ is restricted to $\Omega$ then extended by $0$ outside $\Omega$. Then, let $k \in \zz_{>0}$ be such that $kN_j$ is bigger than the diameter of $\Omega \subset \zz$, then $\wt{y}(m) = \sum_{n \in \Omega} 2y_{kN_j}(m-n) y(n)$ is an element of $B^{Y_j}_1$. Thence $\dl{1} Y_j \geq 1$, and as the other inequality is automatic, $\dl{1} Y_j = 1$.
\par We claim that $Y_\infty = \cap Y_j = \{0\}$. If this were false, then a non-trivial element $y \in Y_\infty$ would have the property that
\[
\forall i \in \zz, \forall n \in \zz, -y(i) = \somme{0\neq k \in \zz}{} y(i+kn).
\]
To get a contradiction, take the limit when $n \to \infty$ and show that it is equal to $0$. First we normalize $y$ so that it is of norm $1$ and suppose that $|y(i)|> \delta$ for some $i$. As an absolutely convergent sequence, $y$ should be concentrated on some set: there exists $n_\delta$ such that $\nr{y}_{\ell^1(\Omega_{n_\delta})} \geq 1- \delta/2$. However when $n > 2n_\delta+1$
\[
|y(i)| = \left| \somme{0\neq k \in \zz}{} y(i+kn) \right| \leq \delta/2,
\]
which is a contradiction.
Thus $\dl{1} Y_\infty = 0$ whereas \mbox{$\limm{n \to \infty} \dl{1} Y_n =1 $}.
\end{ex}
\par Such an possibility is fortunately confined to $\ell^1$; more generally the above construction can be described as follows. Let $\Gamma' \subset \Gamma$ be a subgroup of finite index. Let $\pi:Y \to W$ be a $\Gamma'$-invariant linear map from $Y \subset \ell^p(\Gamma;V)$ to a finite dimensional vector space $W$. Then the existence of such a map implies the existence of $\dim W$ elements of  $\ell^{p'}(\Gamma;V^*)$ which are invariant by $\Gamma'$. This is impossible if $p' \neq \infty$, as such elements would not be decreasing at infinity. For such spaces to exists, $p'$ must be $\infty$.

\section{P1 and Ornstein-Weiss' Lemma}\label{slow}

\hspace*{0.75\parindent} The aim of this section is to show independence (P1) on the choice F{\o}lner sequence. This will be achieved by extending Ornstein-Weiss' lemma to meet our needs.
\begin{teo}\label{low}
Let $\Gamma$ be a discrete amenable group. Let $a: \rr_{\geq 0} \times \jo{P}_{finite}(\Gamma) \to \rr_{\geq 0}$ be a function such that, $\forall \Omega, \Omega' \subset \Gamma$ are finite and $\forall \eps \in \rr_{>0}$
\[
\begin{array}{lll}
(\mathrm{a}) \; a \textrm{ is $\Gamma$-invariant, \ie }
          & \forall \gamma \in \Gamma,
                   & a(\eps,\gamma \Omega) = a(\eps,\Omega)\\
(\mathrm{b}) \; a \textrm{ is decreasing in $\eps$, \ie}
          & \forall \eps'\leq \eps,
                   &  a(\eps',\Omega) \geq a(\eps,\Omega) \\
(\mathrm{c}) \; a \textrm { is $K$-sublinear in $\Omega$, \ie}
          & \exists K \in \rr_{>0},
                   & a(\eps, \Omega) \leq K|\Omega|\\
(\mathrm{d}) \; a \textrm { is $c$-subadditive in $\Omega$, \ie}
          & \exists c \in ]0,1]   ,
                   & a(\eps, \Omega\cup \Omega') \leq a(c\eps, \Omega) + a(c\eps, \Omega')\\
\end{array}
\]
then, for any F{\o}lner sequence $\{\Omega_i\}$,
\[
\limm{\eps \to 0} \lims{i \to \infty} \frac{a(\eps,\Omega_i)}{|\Omega_i|} = \limm{\eps \to 0} \limi{i \to \infty} \frac{a(\eps,\Omega_i)}{|\Omega_i|}.
\]
Furthermore, these limits are independent of the chosen sequence $\{\Omega_i\}$.
\end{teo}
\indent Remark that the $K$-sublinear hypothesis (c) is equivalent to another statement. Indeed, using $c$-subadditivity (d), $\Gamma$-invariance (a) and monotonicity in $\eps$ (b), for all $\Omega$, $a(\eps, \Omega) \leq a(c^{|\Omega|}\eps,e_\Gamma) |\Omega|$ where $e_\Gamma \in \Gamma$ is the neutral element. Thus if (a), (b) and (d) hold, then (c) $\ssi \limm{\eps \to 0} a(\eps,e_\Gamma) < \infty$.
\par This understood, the previous theorem is a generalization of the Ornstein-Weiss lemma. Indeed, the assumptions of the latter, are that $a(\eps, \Omega) = a(\Omega)$ is is sub-additive: then monotonicity (b) always hold, being $K$-sublinear (c) is automatic (see above), and $c$-subadditivity (d) is equivalent to usual subadditivity ($c=1$).
\par The proof of theorem \ref{low} being quite technical, let us first show why P1 is a consequence of this theorem.
\begin{cor}\label{folind}
$\dlp$ is independent of the choice of F{\o}lner sequence.
\end{cor}
\begin{proof}
It suffices to prove that for any $Y \subset \ell^\infty(\Gamma;V)$ a $\Gamma$-invariant set, theorem \ref{low} can be invoked, where $a(\eps,\Omega) = \wdm_\eps (B^{Y,p}_r,\ev_{\ell^p(\Omega)})$. Here is why:
\begin{enumerate}\renewcommand{\labelenumi}{{\normalfont (\alph{enumi})}}
\item By $\Gamma$-invariance of $Y$.
\item As $\wdm_\eps$ is decreasing in $\eps$ (\cf proposition \ref{propbase}.a).
\item Using proposition \ref{propbase1}.b $\wdm_\eps (B^{\ell^p(\Omega)}_r, \ev_{\ell^p(\Omega)}) \leq |\Omega|\dim V$. Then proposition \ref{propbase}.a allows to conclude as $(B^{Y,p}_r, \tau^*, \ev_{\ell^p(\Omega)})$ can be sent without reducing distance by a continuous map to $(B^{\ell^p(\Omega)}_r, \ev_{\ell^p(\Omega)})$
\item Let $\pi_i : \ell^p(\Gamma,V) \to \ell^p(\Omega_i,V)$, then
\[
\pi_1 \times \pi_2 : (B^{Y,p}_r,\ev_{\ell^p(\Omega)}) \to (B^{Y,p}_r,\ev_{\ell^p(\Omega_1)}) \times_p (B^{Y,p}_r,\ev_{\ell^p(\Omega_2)})
\]
is a linear map that does not reduce distances. Applying proposition \ref{propbase}.a and \ref{propbase}.c, yields
\[
\wdm_{2^{1/p}\eps} (B^{Y,p}_r,\ev_{\ell^p(\Omega)}) \leq \wdm_{\eps}(B^{Y,p}_r,\ev_{\ell^p(\Omega_1)}) + \wdm_{\eps} (B^{Y,p}_r,\ev_{\ell^p(\Omega_2)}).
\]
Thence, we conclude that $a(\eps,\Omega)$ is $2^{-1/p}$-subadditive. \hfill \qedhere
\end{enumerate}
\end{proof}
The following notations and definitions will be required in our arguments. The original proof of the Ornstein-Weiss lemma can be found in \cite{OW}. The proof that can be better adapted to our case is however that of \cite[{\S}1.3.1]{Gro} (also explained in \cite{Kri}).
\begin{defi}\label{mutil}
Let $\Gamma$ be a group, let $F \subset \Gamma$ be such that $e_\Gamma \in F$ then  define respectively the outer $F$-boundary, the inner $F$-boundary, the $F$-boundary, the $F$-interior and the $F$-closure of $\Omega$ to be
\[
\begin{array}{lll}
\del_F^+ \Omega
   &= \{\gamma \notin \Omega | \gamma F \cap \Omega \neq \vide \textrm{ and } \gamma F \cap \Omega^c \neq \vide\}
   &= F^{-1} \Omega \cap \Omega^c \\
\del_F^- \Omega
   &= \{\gamma \in \Omega | \gamma F \cap \Omega \neq \vide \textrm{ and } \gamma F \cap \Omega^c \neq \vide\}
   &= F^{-1} \Omega^c \cap \Omega \\
\del_F   \Omega
   &= \{\gamma \in \Gamma | \gamma F \cap \Omega \neq \vide \textrm{and } \gamma F \cap \Omega^c \neq \vide\}
   &= \del_F^+ \Omega \cup \del_F^- \Omega \\
\IN_F   \Omega
   &= \{\gamma \in \Gamma | \gamma F \subset \Omega \}
   &= \Omega \setminus \del^-_F   \Omega \\
\FE_F   \Omega
   &= \{\gamma \in \Gamma | \gamma F \cap \Omega \neq \vide \}
   &= \Omega \cup \del_F^+ \Omega. \\
\end{array}
\]
Moreover, let $|\cdot|$ denote a measure on $\Gamma$. The relative amenability function will be defined as $\alpha(\Omega;F) = \frac{|\del_F\Omega|}{|\Omega|}$, given that these numbers are finite.
\end{defi}
\par Before we move on to technical results, observe that the F{\o}lner condition implies that $\alpha(\Omega_i;F) \to 0$ for any finite set $F$ and any F{\o}lner sequence $\{\Omega_i\}$. Another useful property is that if $F'\subset F$, then $\alpha(\Omega;F') \leq \alpha(\Omega;F)$ since $\del_{F'}\Omega \subset \del_F \Omega$. We start by showing covering properties of big sets by smaller sets.
\begin{defi}\label{dqt}
Let $\eps \in]0,1[$. Subsets $F_i$ of finite measure of $\Gamma$ will be said $\eps$-disjoint if there exists $F'_i \subset F_i$ which are disjoint and such that $|F_i'| \geq (1-\eps)|F_i|$ and $\cup F_i' = \cup F_i$.
\par A subset of finite measure $\Omega$ will be said to admit an $\eps$-quasi-tiling by the subsets $F_i$ if
\begin{enumerate}\renewcommand{\labelenumi}{{\normalfont (\alph{enumi})}}
\item $F_i \subset \Omega$,
\item the $F_i$ are $\eps$-disjoint,
\end{enumerate}
\end{defi}
Here is a first lemma which studies the proportion of a set $\Omega$ covered by an $\eps$-quasi-tiling of translates of another set $F$.
\begin{lem}\label{epac}
Let $\Gamma$ be a discrete group endowed with the counting measure, denoted by $|\cdot|$. Let $\Omega\subset \Gamma$ and $e_\Gamma \in F \subset \Gamma$ both finite sets and such that $\alpha(\Omega;F) < 1$. Let $\{\gamma_i\}_{1\leq i \leq k}$ be a maximal sequence of elements of $\Gamma$ such that the $\gamma_i F$ form an $\eps$-quasi-tiling of $\Omega$. Let $U_F^i = \union{j=1}{i} \gamma_j F$, then
\[
\fr{|U_F^k|}{|\Omega|} \geq \eps (1-\alpha(\Omega;F)).
\]
\end{lem}
\begin{proof}
(This proof corresponds to the first part of the proof of the Ornstein-Weiss lemma in \cite[{\S}1.3.1]{Gro}.) We shall use this general fact:
\[
\begin{array}{rl}
\bint{\Gamma}{} |G_1 \cap \gamma G_2 | \dd\mu(\gamma)
     &= \bint{\Gamma}{} \bint{\Gamma}{} \un_{G_1 \cap \gamma G_2}(\gamma') \dd \mu(\gamma') \dd \mu(\gamma) \\
     &= \bint{\Gamma}{} \bint{\Gamma}{} \un_{G_1}(\gamma') \un_{G_2}(\gamma'\gamma^{-1}) \dd \mu(\gamma) \dd \mu(\gamma')\\
     &= \bint{\Gamma}{} \un_{G_1}(\gamma') \left( \bint{\Gamma}{} \un_{G_2}(\gamma'\gamma^{-1}) \dd \mu(\gamma) \right) \dd \mu(\gamma')\\
     &= \bint{\Gamma}{} \un_{G_1}(\gamma') |G_2| \dd\mu(\gamma')\\
     &= |G_1| |G_2|.
\end{array}
\]
Thus,
\[
\begin{array}{rl}
|\IN_F\Omega|^{-1} \bint{\IN_F\Omega}{} |U_F^k \cap \gamma F| \dd \mu (\gamma)
    & \leq |\IN_F\Omega|^{-1} \bint{\Gamma}{} |U_F^k \cap \gamma F| \dd \mu (\gamma) \\
    & \leq (1-\alpha(\Omega;F))^{-1} |\Omega|^{-1} |U^k_F| |F|.
\end{array}
\]
Clearly, $|U_F^{i-1} \cap \gamma_i F| \leq \eps |F|$, as the $ \gamma_i F $ are $\eps$-disjoint. On the other hand, maximality of $k$ implies that $\forall \gamma \in \IN_F \Omega, |U_F^k \cap \gamma F| \geq \eps |F|$. We then observe that
\[
|\IN_F\Omega|^{-1} \bint{\IN_F\Omega}{} |U_F^k \cap \gamma F| \dd \mu (\gamma) \geq \eps |F|.
\]
Consequently, $\eps (1-\alpha(\Omega;F)) \leq |U^k_F|/|\Omega|$.
\end{proof}
Note that the quasi-tiling can be empty if $\alpha(\Omega;F) = 1$. More precisely, the proof actually works for $\alpha^-(\Omega;F) := \frac{|\del_F^-\Omega|}{|\Omega|}$ instead of $\alpha$. It has the advantage that $\IN_F \Omega \neq \vide$ implies that $\alpha^-(\Omega;F) < 1$ and the quasi-tiling is non-empty. 
In any case, in the upcoming applications, $F$ will always be contained in $\Omega$. The three following lemmas are technical ingredients which will be used in the proof of the generalisation of the Ornstein-Weiss lemma.
\begin{lem}\label{lowtek1}
Let $\Omega' \subset \Omega \subset \Gamma$ and $F \subset \Gamma$ be finite. Suppose that there exists $\eps$ such that \mbox{$|\Omega \setminus \Omega'| \geq \eps |\Omega| $}, then
\[
\alpha(\Omega \setminus \Omega';F) \leq \frac{\alpha(\Omega';F)+\alpha(\Omega;F)}{\eps}.
\]
\end{lem}
\begin{proof}
Since $|\del_F (\Omega \setminus \Omega')| \leq |\del_F \Omega | + |\del_F \Omega'| = \alpha(\Omega;F) |\Omega| + \alpha(\Omega';F) |\Omega'|$, and that \linebreak \mbox{$|\Omega \setminus \Omega'| \geq \eps|\Omega|\geq \eps |\Omega'|$}, a substitution yields
\[
\alpha(\Omega \setminus \Omega';F) = \frac{|\del_F (\Omega \setminus \Omega')|}{|\Omega \setminus \Omega'|} \leq \frac{\alpha(\Omega;F)|\Omega|}{\eps|\Omega|} + \frac{\alpha(\Omega';F) |\Omega'|}{\eps|\Omega'|}.\qedhere
\]
\end{proof}
\begin{lem}\label{lowtek2}
  Let $F \subset \Gamma$ be finite, and let $\{D_i\}_{1\leq i \leq n}$ be an $\eps$-disjoint family of subsets. Then
\[
\alpha(\cup D_i;F) \leq \frac{\max(\alpha(D_i;F))}{1-\eps}
\]
\end{lem}
\begin{proof}
Since $\del_F (\cup D_i) \subset \cup \del_F D_i$, we obtain that
\[
\abs{\del_F (\cup D_i)} \leq \sum \abs{\del_F D_i} \leq \sum \alpha(D_i;F) \abs{D_i} \leq \max(\alpha(D_i;F)) \sum \abs{D_i}
\]
However $(1-\eps)\sum \abs{D_i} \leq \abs{\cup D_i}$ as they are $\eps$-disjoint. Thus
\[
\alpha(\cup D_i;F) = \frac{\del_F (\cup D_i)}{\abs{\cup D_i}} \leq \frac{\max(\alpha(D_i;F))}{1-\eps}. \qedhere
\]
\end{proof}
The last lemma is an adaptation of a useful property of $\zz$ to general amenable group. Consider the typical F{\o}lner sequence for $\zz$, $I_i=[-i,i]$. Then any sufficiently big interval in this family is covered (except for small bits) by translates of some $I_i$.
\begin{lem}\label{paciter}
Let $\{F_i\}$ be a F{\o}lner sequence, let $\delta \in ]0,1/2[$. Then there exists a subsequence (which depends on $\delta$) $\{F_{n_i}\}$, an integer $N(\delta)$, and a sequence of integers $\{k_i\}_{1 \leq i \leq N}$ such that for all set $\Omega$ containing $F_{n_N}$ and satisfying $\alpha(\Omega, F_{n_N}) \leq 2\delta^{2N}$ there exists a  family $\jo{G}$ of $\delta$-disjoint sets such that $|\union{F \in \jo{G}}{} F | \geq (1-\delta)|\Omega|$ and $\jo{G}$ consists in $k_i$ translates of the sets $F_{n_i}$
\end{lem}
\begin{proof}
(This part argument of the argument is briefly presented in \cite[{\S}1.3.1]{Gro}, more details are found in \cite{Kri}; the original result can be found in \cite{OW}.)
In order to better show how the constants enter the proof, we denote $\eps_1= \delta^{2N}$, $\eps_2 = 2\delta^{2N}$ and $\rho= \delta$. First, $\forall \eps_1 \in ]0,1[$, it is possible to refine the sequence $\{F_i\}$ to have
\[
\alpha(F_{i+1},F_i) \leq \eps_1.
\]
Now, let $\Omega^{(1)}=\Omega$ so that $\alpha(\Omega^{(1)}, F_n) \leq \eps_2$, where $n$ will be determined later on. We will cover $\Omega^{(1)}$ to a proportion of $1-\delta$ by almost disjoint translates of the $F_i$, where $1\leq i \leq n$, in $n$ steps (or less). For any $\rho \in ]0,\frac{1}{2}[$, lemma \ref{epac} gives a $\rho$-quasi-tiling of $\Omega^{(1)}$ by $k_n$ translates of $F_n$ such that $| U^{k_n}_{F_n}| \geq \rho (1-\eps_2)|\Omega^{(1)}|$. Let  $\Omega^{(2)} = \Omega^{(1)} \setminus U^{k_n}_{F_n}$, then $|\Omega^{(2)}| \leq (1-\rho+\eps_2\rho)|\Omega^{(1)}|$.
\par If $|\Omega^{(2)}| \leq \delta |\Omega^{(1)}|$ the goal is achieved and there is no need to continue. Otherwise, lemma \ref{lowtek1} then lemma \ref{lowtek2} shows that
\[
\alpha(\Omega^{(2)}, F_{n-1}) \leq \frac{1}{\delta}(\eps_1 + \alpha(U^{k_n}_{F_n};F_{n-1}) ) \leq \frac{1}{\delta}(\eps_1 +\frac{\eps_1}{1-\rho}) \leq 3  \frac{\eps_1}{\delta}.
\]
It is now possible to recover $\Omega^{(2)}$ by a $\rho$-quasi-tiling of $k_{n-1}$ translates of $F_{n-1}$ in such a way that $| U^{k_{n-1}}_{F_{n-1}}| \geq \rho (1- 3 \frac{\eps_1}{\delta})|\Omega^{(2)}|$. We now have a set $\Omega^{(3)}$ such that
\[
|\Omega^{(3)}| \leq (1- \rho- 3 \rho \frac{\eps_1}{\delta} )|\Omega^{(2)}| \leq (1-\rho+\eps_2\rho)(1- \rho+ \rho \frac{3\eps_1}{\delta} )|\Omega^{(1)}|
\]
We will now take $\eps_2 = 2 \eps_1$. Proceeding by induction, as long as $|\Omega^{(i-1)}|\geq \eps |\Omega^{(1)}|$, the set $\Omega^{(i)}$ (for $1\leq i \leq n$) will have the following properties:
\begin{enumerate}\renewcommand{\labelenumi}{{\normalfont \arabic{enumi}.}}
\item $\alpha(\Omega^{(i)},F_{n-i+1}) \leq (1+i)\eps_1/\delta^{i-1}$
\item $U_{F_{n-i+1}}^{k_{n-i+1}}$ is a $\rho$-quasi-tiling of $\Omega^{(i)}$ by translates of $F_{n-i+1}$
\item If $\Omega^{(i+1)} = \Omega^{(i)} \setminus U_{F_{n-i+1}}^{k_{n-i+1}}$ then $|\Omega^{(i+1)}| \leq |\Omega^{(1)}|\produ{j=1}{i} (1-\rho (1-(1+j)\eps_1/\delta^{j-1}) )$
\end{enumerate}
Since it is not possible to hope that this process terminates before $i=n$, it remains to be checked that if $n$ is big enough, we still get a quasi-tiling that covers $(1-\delta)|\Omega^{(1)}|$ elements. To achieve this, observe that the product in the third property above can be bounded if $i=n$ by
\[
\produ{j=1}{n} (1-\rho (1-(1+j)\eps_1/\delta^{j-1}) ) \leq (1-\rho (1-(1+n)\eps_1/\delta^{n-1}) )^n.
\]
For $\eps_1 = \delta^{2n}$, the right-hand term tends to $0$ when $n$ tends to $\infty$. Thus, $\exists N(\delta, \rho)$ such that if $\eps_1 = \delta^{2N}$ translates of $F_j$ (where $1 \leq j \leq N$) form a $\rho$-quasi-tiling of any set $\Omega^{(1)}$ such that $\alpha(\Omega^{(1)};F_N) \leq \delta^{2N}$.
\par We substitute as promised $\rho=\delta$ to have: for any fixed $\delta$, we choose  a subsequence whose members satisfy $\alpha(F_{n_{i+1}},F_{n_i}) \leq \delta^{2N}$ where $N$ is such that
\[
(1-\delta (1-(1+N)\delta^{N+1}))^N < \delta.
\]
Then successive applications of lemma \ref{epac} give the required translates of $F_{n_i}$.
\end{proof}
We are now ready to prove the main result of this section. At a first reading, this proof might be easier to understand with $\Gamma = \zz$ in  mind (taking $\Omega_n = [-n,n] \cap \zz$).
\begin{proof}[Proof of theorem \ref{low}]
Let us first introduce some notations for the functions given by pointwise convergence and their limits. Let$\{\Omega_i^{+,\eps}\}$ and $\{\Omega_i^{-,\eps}\}$ be subsequences of $\{ \Omega_i \}$ such that
\[
\limm{i \to \infty} \frac{a(\eps,\Omega_i^{+,\eps})}{|\Omega_i^{+,\eps}|} = \lims{i \to \infty} \frac{a(\eps,\Omega_i)}{|\Omega_i|} \quad \textrm{ and } \quad \limm{i \to \infty} \frac{a(\eps,\Omega_i^{-,\eps})}{|\Omega_i^{-,\eps}|} = \limi{i \to \infty} \frac{a(\eps,\Omega_i)}{|\Omega_i|}.
\]
Using (c), these limits are respectively real numbers  $l^+(\eps)$ and $l^-(\eps)$ belonging to the interval $[0, K]$. Furthermore, let
\[
 l^+ := \limm{\eps \to 0} l^+(\eps) \quad \textrm{ and } \quad l^- := \limm{\eps \to 0} l^-(\eps).
\]
Trivially, $l^+(\eps) \geq l^-(\eps)$, but nothing forces $l^\pm(0)=l^\pm$ (in general, equality is not expected). If we try to use the usual argument directly, a problem arises due to the \linebreak $c$-subadditivity. Indeed, taking a sequence which converges to $l^+(\eps)$ and decomposing it using another sequence which converges to $l^-(\eps)$ by subadditivity will fail. A factor of $c$ will appear in front of the $\eps$ (see$(d)$), and this would force to pass from the sequence $\Omega^{-,\eps}_i$ to $\Omega^{-,c\eps}_i$ at each step. Diagonal arguments settle this problem. Let
\[
b_i(\eps) = \frac{a(\eps,\Omega_i^{-,1/i})}{|\Omega_i^{-,1/i}|}.
\]
This is a sequence of bounded decreasing functions defined for $\eps \in [0,1]$ with value in $[0,K]$. Using (one of) Helly's theorem (\cf \cite[{\S}36.5 theorem 5, p.372]{KF}), there exists a subsequence $n_i$ which possesses a limit at each point. We briefly recall how this subsequence is obtained. First, a sequence $\{r_k\}_{k \geq 1}$ of dense rational number in $[0,1]$ is taken. Since the $b_i(\eps)$ are bounded, let $n_i^{(j)}$ be the subsequence which converges at $r_k$ for $1 \leq k\leq j$. The diagonal sequence $n_i = n_i^{(i)}$ converges at each $r_k$, and since the functions $b_i(\eps)$ are decreasing, the function $l^H(\eps) = \limm{i \to \infty} b_{n_i}(\eps)$ which is {\it a priori} only defined for the $r_k$ is also decreasing. It remains to be checked that $l^H(\eps)$ extended at all the points of $[0,1]$ by approximating by a sequence of increasing $r_k$ is the actual limit of the subsequence $n_i$ (see the above reference for details). Let us show that $\limm{\eps \to 0} l^H(\eps) = l^-$. This follows from
\[
\begin{array}{lll}
\forall \delta > 0, \exists N_1(\delta) \textrm{ such that } N_1(\delta) < i
     &\imp |b_{n_i}(\frac{1}{n_i}) - l^H(\frac{1}{n_i})| < \delta;\\
\forall \delta > 0, \exists N_2(\delta) \textrm{ such that } N_2(\delta) < i
     &\imp |b_{n_i}(\frac{1}{n_i}) - l^-(\frac{1}{n_i})| < \delta;\\
\forall \delta > 0, \exists N_3(\delta) \textrm{ such that } N_3(\delta) < i
     &\imp |l^-(\frac{1}{n_i})-l^-| < \delta .\\
l^H(\eps) \textrm{ is decreasing in } \eps & \imp \limm{\eps \to 0} l^H(\eps) = \limm{i \to \infty} l^H(\frac{1}{n_i})
\end{array}
\]
These four assertions are respectively consequences of the definition of $l^H$, the choice of $\Omega^{-,1/i}_i$, the definition of $l^-$, and the fact that a limit that exists (thanks to monotonicity) is achieved by any sequence. We shall now show that
\[
\forall \delta > 0, l^+(\eps) \leq \limm{\eps' \to 0} l^H(\eps') + \delta = l^- + \delta.
\]
The argument is in essence the same as for subadditive sequences of real numbers: lemma \ref{paciter} plays the role of the decomposition $n = kn' + r$ and $c$-subadditivity (d) forces $\eps \to 0$.
\par Let $\delta \in ]0,\frac{1}{2}[$. Denote by $F_i = \Omega_{n_i}^{-,1/n_i}$. It is possible to refine this sequence so that
\[
a(\eps,F_i)/|F_i| \leq l^H(\eps)+\delta.
\]
Applying lemma \ref{paciter} gives an $\eps$-quasi-tiling (which does not cover a set of proportion $\delta$) of any sufficiently big set by translates of the $F_i$. Since $\{\Omega^{+,\eps}_i\}$ is also a F{\o}lner sequence, for $i$ big enough, lemma \ref{paciter} applies to each element. Take $\Omega= \Omega_i^{+,\eps}$, denote $\gamma_{F_j;m} F_j$ the $k_j$ translates of $F_j$ obtained ($m=1, \ldots, k_j$), and let $i_0$ be such that $|\Omega^{(i_0)}|\leq \delta |\Omega^{(1)}|$. Thanks to repeated use of $c$-subadditivity (d), we have that
\[
\begin{array}{rll}
a(\eps, \Omega^{(1)})
   & \leq & \Somme{m=1}{k_n} a(c^m\eps, \gamma_{F_n;m} F_n) + a(c^{k_n}\eps, \Omega^{(2)})\\
   & \leq & \ldots \\
   & \leq & \Somme{i=n-i_0}{n} \Big( \Somme{m=1}{k_i} a(c^{\kappa_i+m}\eps, \gamma_{F_i;m}F_i) \Big) + a(c^{\kappa_{i_0}}\eps,\Omega^{(i_0)}),
\end{array}
\]
where $\kappa_i = \somme{j=n-i}{n}k_n$. Using $\Gamma$-invariance (a), the fact that these functions are decreasing in $\eps$ (b), and the $K$-sublinear property (c), this inequality yields
\[
a(\eps, \Omega^{(1)}) \leq \somme{i=n-i_0}{n} \Big( \somme{m=1}{k_i} a(c^{\kappa_{i_0}}\eps, F_i) \Big) +K|\Omega^{(i_0)}|.
\]
On one hand, $|\Omega^{(i_0)}|\leq \delta |\Omega^{(1)}|$ and $\frac{a(c^{\kappa_{i_0}}\eps, F_i)}{|F_i|}\leq l^H(c^{\kappa_{i_0}})+\delta $. Thence,
\[
\begin{array}{rll}
\fr{a(\eps, \Omega^{(1)})}{|\Omega^{(1)}|}
  & \leq & \Somme{i,m}{} \fr{a(c^{\kappa_{i_0}}\eps, F_i)}{|F_i|} \fr{|\gamma_{F_i;m}F_i|}{|\Omega^{(1)}|} +K\fr{|\Omega^{(i_0)}|}{|\Omega^{(1)}|} \\
  & \leq & (l^H(c^{\kappa_{i_0}})+\delta) \Somme{i,m}{} \fr{|\gamma_{F_i;m}F_i|}{|\Omega^{(1)}|}+ K \delta
\end{array}
\]
On the other hand, the $\{ \gamma_{F_i;m}F_i \}$ are $\delta$-disjoint. Thus
\[
(1-\delta) \sum | \gamma_{F_i;m}F_i| \leq |\cup \gamma_{F_i;m}F_i| \leq |\Omega^{(1)}|.
\]
This shows that
\[
\frac{a(\eps, \Omega^{+,\eps}_j)}{|\Omega^{+,\eps}_j|} \leq (l^H(c^{\kappa_{i_0}})+\delta) \somme{i,l}{} \frac{|\gamma_{F_i;l}F_i|}{|\Omega^{(1)}|}+ K \delta \leq \frac{l^H(c^{\kappa_{i_0}})+\delta}{1-\delta} + K \delta,
\]
For all $\Omega^{+,\eps}_j$ big enough, where $\kappa_{i_0}$ depends on $\Omega^{+,\eps}_j$. Since $l^H(\eps)$ is decreasing and $\limm{\eps \to 0} l^H(\eps) = l^-$, taking the limit when $j$ and $\kappa_{i_0} \to \infty$ is not  a problem:
\[
l^+(\eps) \leq l^- + \delta (K+l^-+1).
\]
We have shown that $l^+=l^-$.
\par To deduce the independence on the choice of sequence, notice that given two F{\o}lner sequences $\{\Omega_i\}$ and $\{\Omega'_i\}$, the sequence $\{\wt{\Omega}_i\} $ whose elements alternate between those of the two former sequences will also possess a limit. The limit obtained with $\{\srl{\Omega}_i\} $ must be equal to the one taken via $\{\Omega_i\}$ or $\{\Omega'_i\}$.
\end{proof}

\appendix
\section{Von Neumann's dimension and $\dl{2}$}  \label{sdl2}

We recall an argument of Gromov (see \cite[{\S}1.12]{Gro}) that relates von Neumann to the semi-axis of ellipsoids and thus showing that $\dl{2}$ is indeed von Neumann dimension. We briefly review the definition of the latter.
\par Let $Y \subset \ell^2(\Gamma; V)$ be a $\Gamma$-invariant linear subspace, $\forall \Omega \subset \Gamma$ we define the operator $R_\Omega: Y \to \ell^2(\Omega;V)$ by restriction to $\Omega$: $y \mapsto y_{|\Omega}$. Its adjoint $R^*_\Omega: \ell^2(\Omega;V) \to Y$ is the orthogonal projection to $Y$. To see this, write $R_\Omega(y) = y \un_\Omega $ where $\un_\Omega$ is the characteristic function of $\Omega$, then
\[
\gen{R^*_\Omega(x),y} := \gen{x ,R_\Omega y } = \bint{\Gamma}{} x y \un_\Omega = \bint{\Gamma}{} (\un_\Omega x) y.
\]
However this last expression is simply the scalar product of $x$, extended as a function on all of $\Gamma$ by $0$, with $y$. Thus, $R^*_\Omega(x)$ is the projection on $Y$ of the extension of $x$ to $\Gamma$ by $0$. In what follows we will omit this inclusion (extension by $0$) from $ \ell^2(\Omega;V)$ to $ \ell^2(\Omega';V)$ when $\Omega \subset \Omega' $. Dependence on $\Omega$ of $R^*_\Omega$ will not be written. A crucial remark is that the invariance of $Y$ by $\Gamma$ implies that, for $\Omega, \Omega' \subset \Gamma $ finite subsets,
\[
\fr{\tr R_\Omega R^*}{\tr R_{\Omega'} R^*} = \fr{|\Omega|}{|\Omega'|}.
\]
A possible definition of von Neumann dimension (see \cite{Lu} or \cite[{\S}1]{PP}) is
\[
\dl{2} (Y:\Gamma) := |\Omega|^{-1} \tr R_\Omega R^*
\]
for a $\Omega \subset \Gamma$. This quantity is actually independent of the chosen set. The aim of this section is to retrieve this quantity as the $\wdm$ of a certain object.
\begin{teo}\label{vpVN} (\cf \cite[1.12A]{Gro})
Let $\Omega_i \subset \Gamma$ be a F{\o}lner sequence, let $n_i[a,b]$ be the number of eigenvalues of the operator $R_{\Omega_i} R^*$ (defined relative to $Y$) contained in the interval $[a,b]$. If $0<a\leq b<1$, then
\[
\limm{i \to \infty} \fr{n_i [a,b]}{|\Omega_i|} =0
\]
\end{teo}
\begin{proof}
(The proof is with minor differences in notation that of \cite{Gro}.) Since $R_{\Omega}$ and $R^*$ are both projections (in $\ell^2$), the eigenvalues of $R_{\Omega} R^*$ will be contained in $[0,1]$. The proof proceeds in three steps.
\par {\it First,} let $x \in \ell^2(\Omega;V)$, it will be called an $\eps$-quasimode of eigenvalue $\lambda$ for $R_{\Omega} R^*$ if
\begin{equation}\label{qmd}
\nr{R_{\Omega} R^* x - \lambda x}_{\ell^2} \leq \eps \nr{x}_{\ell^2}.
\end{equation}
If $x$ is such an element, and if its restriction outside $\Omega$ is small, more precisely
\begin{equation}
 \label{ppa}
\nr{R^* x_{| \Gamma \setminus \Omega}}_{\ell^2} = \nr{R^* x - R_\Omega R^*x}_{\ell^2} \leq \delta \nr{x}_{\ell^2},
\end{equation}
then $\lambda(1-\lambda) \leq 2 \eps + \delta$. Indeed, using \eqref{qmd} in \eqref{ppa} yields that
\[
\nr{R^* x -\lambda x}_{\ell^2} \leq (\delta + \eps) \nr{x}_{\ell^2}.
\]
$R^*$ is a projection, $R^*R^* = R^*$ and $\nr{R^*} = 1$, thence
\[
\begin{array}{rll}
(1-\lambda) \nr{R^* x}_{\ell^2}
  &= \nr{R^* x- R^* \lambda x}_{\ell^2}
  &=\nr{R^*(R^*x-\lambda x)}_{\ell^2} \\
  & \leq \nr{R^*x-\lambda x}_{\ell^2}
  & \leq (\delta + \eps) \nr{x}_{\ell^2},
\end{array}
\]
as the eigenvalues of $R_\Omega R^*$ are all contained in $[0,1]$, $|1-\lambda| = 1-\lambda$. Moreover the restriction to $\Omega$ can only reduce the norm, $ \nr{R_\Omega R^* x}_{\ell^2} \leq  \nr{R^* x}_{\ell^2}$. Using \eqref{qmd} anew gives,
\[
(1-\lambda) \lambda \nr{x}_{\ell^2} \leq (1-\lambda) \nr{R_\Omega R^* x}_{\ell^2} + (1-\lambda)\eps \nr{x}_{\ell^2} \leq (\delta + (2-\lambda) \eps ) \nr{x}_{\ell^2}
\]
{\it Second,} 
we will look at smaller set inside $\Omega$, see definition \ref{mutil} for notations. Let $F_\rho$ be a set of cardinality $\rho$.
The next argument will consist in showing that most of $x \in \ell^2(\IN_{F_\rho} \Omega; V)$ have a small projection to $\Gamma \setminus \Omega$. That is, let
\[
S_\rho  = R_{\Gamma \setminus \Omega} R^* = (1-R_{\Omega})R^*: \ell^2(\IN_{F_\rho} \Omega;V) \to \ell^2(\Gamma \setminus \Omega;V),
\]
then, for some $F_\rho$, $\tr S^*_\rho S_\rho \leq \beta(\rho) \dim V |\IN_{F_\rho} \Omega|$ where $\beta(\rho)$ tends to $0$ when $\rho \to \infty$. The dependence on $\rho$ does not only come from the domain of definition: the operator $S^*_\rho S_\rho$ is
\[
\begin{array}{rl}
S^*_\rho S_\rho &= R_{\IN_{F_\rho} \Omega} (1-R^*)(1-R_\Omega)R^* \\
                &= (R_{\IN_{F_\rho} \Omega} - R_{\IN_{F_\rho} \Omega} R^*)(R^*-R_\Omega R^*) \\
                &= (R_{\IN_{F_\rho} \Omega} R^* - R_{\IN_{F_\rho} \Omega} R^* R^* - R_{\IN_{F_\rho} \Omega} R_\Omega R^* + R_{\IN_{F_\rho} \Omega} R^* R_\Omega R^*)\\
                &= R_{\IN_{F_\rho} \Omega} R^* (R_\Omega R^* -1).
\end{array}
\]
\par Any Dirac mass $x_\gamma$ with support at a point $\gamma$ satisfies $\nr{R^* x_\gamma}_{\ell^2} \leq 1$ ($R^*$ is a projection). Thus, there exists $F_\rho \ni e_\Gamma$, such that $\pnr{(1-R_{\gamma F_\rho}) R^*x_\gamma}_{\ell^2} \leq \beta(\rho)$ where $|F_\rho|=\rho$ and $\limm{\rho \to \infty} \beta(\rho) =0$. Consequently, for $\gamma \in \IN_{F_\rho} \Omega$, $\nr{S_\rho x_\gamma}_{\ell^2} \leq \beta(\rho)$ since $\Gamma \setminus \Omega$ is contained in the complement of the union of the $\gamma F_\rho$ for $\gamma \in \IN_{F_\rho} \Omega$. Since $\pnr{S_\rho^*}\leq 1$, then $\pnr{S_\rho^* S_\rho x_\gamma}_{\ell^2} \leq \beta(\rho)$. The Dirac masses being an orthonormal basis for $\ell^2(\IN_{F_\rho} \Omega;V)$, it follows that $\tr S^*_\rho S_\rho \beta(\rho) \leq \beta(\rho) \dim V  |\IN_{F_\rho} \Omega| $.
\par {\it Last,} we shall evaluate $n_i[a,b]$ for $a,b \in ]0,1[$ and $b-a=\eps \in ]0,1[$. Let $X_i$ be the space generated by eigenvectors of $R_{\Omega_i}R^*$ whose eigenvalue is in $[a,b]$. Then, $\forall \lambda \in [a,b], \forall x \in X_i$, $x$ is an $\eps$-quasimode of eigenvalue $\lambda$ for $R_{\Omega_i}R^*$. The evaluation of $\dim X_i$ will be done by looking at spaces whose dimension is close. Let \linebreak $X_i^\rho = X_i \cap \ell^2(\IN_{F_\rho} \Omega_i,V)$ be the subspace of elements which vanish on the thickened boundary, $\dim X_i - \dim X_i^\rho \leq \dim V|\del^-_{F_\rho} \Omega_i|$. The amenability of $\Omega_i$ implies that this difference is negligible, $\limm{i \to \infty} (\dim X_i - \dim X_i^\rho)/|\Omega_i| = 0$; it will suffice to evaluate $\dim X_i^\rho$.
\par Unfortunately, neither $X_i^\rho$ nor $X_i$ is {\it a priori} invariant by $S_\rho^* S_\rho$. Let's nevertheless look at the intersection of $X_i^\rho$ with the space generated by eigenvectors of $S_\rho^* S_\rho$ of eigenvalue $\leq \kappa^2$; we will denote this new intersection by $X_i^{\rho,\kappa}$. On this space, $\nr{S_\rho} \leq \kappa$ since
\[
\nr{S_\rho x}_{\ell^2}^2 = \pgen{S_\rho x, S_\rho x} = \pgen{x ,S_\rho^* S_\rho x} \leq \kappa^2 \nr{x}_{\ell^2}^2
\]
Yet again, this space is of dimension close to that of $X_i^\rho$: if $E^{>\kappa^2}$ is the space of eigenvectors of $S_\rho^* S_\rho$ whose eigenvalue is greater than $\kappa^2$, then
\[
\dim X_i^\rho - \dim X_i^{\rho, \kappa} \leq \dim E^{>\kappa^2} \leq \kappa^{-2} \tr S_\rho^* S_\rho \leq s |\Omega_i^{-\rho}| \beta(\rho) / \kappa^2.
\]
In other words,
\[
\forall \kappa >0, \forall \alpha>0, \exists \rho \textrm{ such that } \lims{i\to \infty} \fr{\dim X_i^\rho - \dim X_i^{\rho,\kappa}}{|\Omega_i|} \leq \alpha,
\]
Thus, it remains to evaluate $\dim X_i^{\rho,\kappa}$. To do so, we use the conclusion of the first part for $\lambda = a$, $\delta = \kappa$ and $\eps = b-a$: if $\dim X_i^{\rho,\kappa} >0$, then $a(1-a) \leq 2(b-a)+\kappa$. Consequently, the inequality $b-a < (a-a^2-\kappa)/2 $ implies that $\dim X_i^{\rho,\kappa}=0$. So when $\rho\geq \rho_0(\kappa,\alpha)$ is sufficiently big, $\lims{i \to \infty} \dim X_i^\rho/|\Omega_i|\leq \alpha$. It follows that
\[
\begin{array}{rl}
\lims{i \to \infty} \fr{\dim X_i}{|\Omega_i|}
  &\leq \lims{i \to \infty} \fr{\dim X_i - \dim X_i^\rho}{|\Omega_i|} + \lims{i \to \infty} \fr{\dim X_i^\rho}{|\Omega_i|} \leq \alpha. \\
\end{array}
\]
But $\alpha \to 0 $ as $\rho \to \infty$. This proves the theorem for intervals $[a,b]$ satisfying $b-a < a(1-a)/2$, as $\kappa$ can also be made arbitrarily small for $\rho$ big enough. The conclusion is obtained by noticing that any interval strictly contained in $[0,1]$ can be covered by intervals of this type.
\end{proof}
\par This property enables us interpret von Neumann dimension as a $\wdm$ for a set with a chosen pseudo-metric.
\begin{cor}\label{wdmVN} (\cf \cite[corollary 1.12.2]{Gro})
Let $Y \subset \ell^2(\Gamma;V)$ be an invariant subspace, let $B^Y_1 = Y \cap B^{\ell^2(\Gamma;V)}_1$ the intersection of the unit ball with $Y$. Then, for a given F{\o}lner sequence $\Omega_i \subset \Gamma$,
\[
\forall \eps \in ]0,1[, \qquad \limm{i \to \infty} \frac{1}{|\Omega_i|} \wdm_\eps (R_{\Omega_i} B^Y_1 , \ell^2 )=  \dl{2} Y
\]
\end{cor}
\begin{proof} (We give the argument of \cite{Gro} in detail.)
To get this result $R_\Omega B^Y_1$ must be seen as an ellipsoid whose semi-axes are related to the eigenvalues of $R_\Omega R^*$. Remark that $B^Y_1= R^* B^{\ell^2(\Gamma;V)}_1$. Then, an ellipsoid can be defined as the image of a ball by an self-adjoint operator, say $A$; the semi-axis of this ellipsoid are in correspondence with the eigenvalue of $A$. It might be worth recalling how this relates to the usual definition of an ellipsoid $E$ (as the set $\{ y  | \gen{y,Py}\leq 1 \}$ for a positive definite operator $P$). The semi-axes of $E$ are of the form $\lambda_i(P)^{-1/2}$ for $\lambda_i(P)$ an eigenvalue of $P$. Indeed let $B^V$ be a ball in a vector space $V$, and let $A:V \to V$ be self-adjoint. Restricting to $V' = \img A = \ker A ^\perp \subset V$, it must be shown that for $x \in V'$ such that $\gen{x,x} \leq 1$, there exists $P:V' \to V'$ positive definite such that $\gen{Ax,PAx} \leq 1$. Taking $P = A^{-2}$ yields the conclusion: $A^{-2}$ is a positive definite operator on $V'$ whose eigenvalues are $\lambda_i(A)^{-2}$. Thus $A B^V$ is an ellipsoid with semi-axis $\lambda_i(P)^{-1/2} = \lambda_i(A)$.
\par In our present context, $R_\Omega R^*$ is self-adjoint, thus $R_\Omega R^* B^{\ell^2(\Gamma;V)}_1 = R_\Omega B^Y_1$ is an ellipsoid whose semi-axis are the eigenvalues of $R_\Omega R^*$. This ellipsoid contains isometrically the ball obtained by ignoring the semi-axis of length $<\eps$ and replacing the remaining ones by semi-axis of length $\eps$. Thus $\wdm_\eps (R_{\Omega_i} B^Y_1, \ell^2) \geq n_i[\eps,1]$. On the other hand, $\wdm_\eps (R_{\Omega_i} B^Y_1, \ell^2 ) \leq n_i[\eps/2,1]$, as the continuous map obtained by projecting on the sub-ellipsoid formed by the semi-axis of length $>\eps/2$ indicates. When $i \to \infty$, the eigenvalues of $R_{\Omega_i} R^*$ tend to $0$ or $1$. In particular, when $i \to \infty$ the inequality
\[
\frac{1}{|\Omega_i|} n_i[\eps,1] \leq  \frac{1}{|\Omega_i|} \wdm_\eps (R_{\Omega_i} B , \ell^2 ) \leq \frac{1}{|\Omega_i|} n_i[\eps/2,1]
\]
shows that $\limm{i \to \infty} \frac{1}{|\Omega_i|} \wdm_\eps (R_{\Omega_i} B^Y_1 , \ell^2 ) =  \dl{2} Y$, since  $n_i[a,1] \to \tr R_{\Omega_i} R^*$.
\end{proof}
This corollary can be expressed in terms of $\ell^p$ dimension. Indeed, let \linebreak $B^Y_1 = Y \cap B_1^{\ell^2(\Gamma;V)}$ be endowed with the pseudo-metric of evaluation at $e \in \Gamma$: $\ev(x,y) = \nr{x(e)-y(e)}_{V}$. Translation of this pseudo-metric by an element of $\gamma$ is the evaluation at $\gamma$. Thus, $\ev_{\ell^2(\Omega)}(x,y) = \nr{x-y}_{\ell^2(\Omega)} = \nr{R_\Omega(x-y)}_{\ell^2}$. The map $R_\Omega: B^Y_1 \to R_\Omega B$ is continuous for the topology of $B^Y_1$ as a subset of $\ell^p$ (with $\tau^*$ or even with the norm topology). 
The fibers are of ``diameter'' $0$ given that $\Omega' \subset \Omega$. Thus, corollary \ref{wdmVN} can be expressed as follows:
\[
\mme_{\ell^2} (B^Y_1,\tau^*,\ev,\{ \Omega_i \}) = \dl{2} Y.
\]
Indeed, $R_{\Omega_i} B^Y_1$ injects isometrically in $(B^Y_1, \ev_{\Omega_i})$ and $(B^Y_1, \ev_{\Omega_i})$ possesses a map to $R_{\Omega_i} B^Y_1$ whose fiber are of ``diameter'' $0$. Thus
\[
\wdm_\eps (B^Y_1, \ev_{\Omega_i}) = \wdm_\eps (R_{\Omega_i} B^Y_1 , \ell^2 ).
\]
 This shows that definition \ref{ddlp} is equivalent when $p=2$ to the von Neumann dimension and this for any F{\o}lner sequence $\{\Omega_i\}$ chosen.
\par It would have been surprising that this were not the case in general. An alteration of the Ornstein-Weiss lemma (see section \ref{slow}) enables to show the independence of the limit on the sequence chosen.

\medskip
Received May 2009; revised September 2009.
\medskip

\end{document}